# Analytical and numerical solutions to ergodic control problems arising in environmental management

**(Running title: ergodic environmental management)**

**Authors** Hidekazu Yoshioka[1, *], Motoh Tsujimura[2], Yuta Yaegashi[3]

[1] Assistant Professor, Graduate School of Natural Science and Technology, Shimane University, Nishikawatsu-cho 1060, Matsue, 690-8504, Japan
[2] Professor, Graduate School of Commerce, Doshisha University, Karasuma-Higashi-iru, Imadegawa-dori, Kamigyo-ku, Kyoto, 602-8580, Japan
[3] Independent Researcher, Dr. of Agr., 10-12-403, Maeda-cho, Niihama, 792-0007, Japan
* Corresponding author
   E-mail: yoshih@life.shimane-u.ac.jp

**Abstract** Environmental management optimizing a long-run objective is an ergodic control problem whose resolution can be achieved by solving an associated non-local Hamilton–Jacobi–Bellman (HJB) equation having an effective Hamiltonian. Focusing on sediment storage management as a modern engineering problem, we formulate, analyze, and compute a new ergodic control problem under discrete observations: a simple but non-trivial mathematical problem. We give optimality and comparison results of the corresponding HJB equation having unique non-smoothness and discontinuity. To numerically compute HJB equations, we propose a new fast-sweep method resorting to neither pseudo-time integration nor vanishing discount. The optimal policy and the effective Hamiltonian are then computed simultaneously. Convergence rate of numerical solutions is computationally analyzed. An advanced robust control counterpart where the dynamics involve uncertainties is also numerically considered.

## 1. Introduction

We consider and mathematically as well as numerically analyze a new ergodic control model based on an integro-differential equation for managing environment under anthropogenic pressure. Environmental management always requires a long-run viewpoint to achieve co-existence with environment anthropogenic pressure. Especially, sustainable management of river environment has been paid much attention recently because of the urgent need to find the way to mitigate the impacts of human threats such as dam and weir constructions [1], environmental pollution [2], and water abstractions [3].

Mathematical approaches have been used for finding the way to set up sustainable environmental management policies. Optimal control theory has been a central mathematical tool for this purpose [4]. Sustainable policies should be robust so that environmental and ecological dynamics, which are often stochastic [5] due to being exposed to uncertain external disturbances [6], are successfully managed. Feedback controls are preferred for managing these dynamics so that they are adaptively controlled to be close to the desired states [7].

Ergodic control is a mathematical concept for adaptively and sustainably controlling stochastic dynamics in a long run [8]. In an ergodic control problem, the target dynamics are optimized so that a time-averaged objective function is minimized or maximized. This concept is suited to many real-world problems, including harvesting of renewable resources [9], chemostat control [10], operation of dam-reservoir systems [11], investment for sustainable economic growth [12], welfare accumulation [13], inventory control [14], mechanical systems management [15], physiological dynamics [16], and hydrodynamic problems [17]. Modeling environmental and ecological dynamics in rivers with the ergodic control, however, has been paid much less attention so far.

In most cases, solving an ergodic control problem is based on a dynamic programming approach that reduces a problem to finding a proper solution to an optimality equation called Hamilton–Jacobi–Bellman (HJB) equation. HJB equations are non-linear and sometimes non-local differential equations that are solvable analytically only in the limited cases [18-22] and must be handled numerically in general.

Therefore, the core of this approach is to find or approximate solutions to the HJB equations. Approaches other than the dynamic programing such as the linear programming [23] and backward stochastic differential equation (SDE) approach [24] are also available. Numerical schemes for degenerate elliptic and/or parabolic equations have been applied to computing HJB equations in ergodic and other stochastic control problems. They include but are not limited to the finite difference [25, 26], finite element [27], and Markov chain approximation schemes [28]. Other numerical schemes like the semi-Lagrangian [29, 30], fitted finite volume [31, 32], and discontinuous Galerkin schemes [33, 34] can also be utilized.

We focus on sediment replenishment as a control problem [35]. The core of the problem is simple; constructing a dam across a river cross-section blocks longitudinal transport of sediment (soil and gravel particles) from the upstream river. This physical phenomenon is called sediment trapping [36]. River environment without sediment supply from the upstream encounters a variety of negative impacts such as severe river erosion and loss of habitats of native species [36]. Sediment replenishment is a mitigation measure against the sediment trapping that transports earth and soils to the target river from another place. There are at least two goals in a sediment replenishment problem: avoiding sediment depletion while minimizing replenishment cost. The central optimization problem of sediment replenishment is thus to find when and how much of sediment should be replenished. These issues have been experimentally investigated in fields and laboratories [37-41]; however, theoretical aspects of the problem have been less explored to the best of the authors' knowledge.

The objective of this paper is to formulate and analyze a sediment replenishment problem as a simple but non-trivial long-run (ergodic) management problem involving non-smooth and jump-driven dynamics, discrete and random observations, and impulse control. This is a new optimization problem that has not been addressed so far despite their importance. We consider that sediment storage follows a non-smooth deteriorating process whose deterioration, decrease of the sediment storage, is driven by deterministic and jump processes [42-44]. The former represents the gradual decrease by base river flows, while the latter represents sudden decrease by floods [45]. The non-smooth property comes from a constraint that the sediment storage should be a non-negative physical variable. This kind of boundedness constraints are often used for modeling and control in environmental management such as storage systems operation [11, 46] and water quality management [47]. The non-negativity constraint is enforced by truncating the drift and jump, leading to unique non-smooth coefficients of the dynamics.

We assume that the driving jump process is a Lévy process having infinite activities [48] to describe generic river flows in a simple way. This is a critical difference between the earlier models where only the compound Poisson noises were considered [43, 44]. Sediment transport is a complex physical process whose complete physical description is still unclear and is often intermittent and fractional [49], motivating us to use the versatile jump processes having infinite activities. Sediment replenishment problems with Lévy processes having infinite activities as flexible models have not been considered so far.

The objective function in our model is a time-averaged value of the sum of a penalty for the sediment depletion and the total cost of sediment replenishment. This penalization induces another non-smooth property into the model. We assume that the opportunities of observation and sediment replenishment arrive only discretely following a Poisson process. This assumption simplifies the mathematical modeling in stochastic control and related problems with discrete and random observations [50-53]. The resulting HJB equation has unique non-smooth and discontinuous coefficients due to the non-smooth properties of the dynamics and objective function. Nevertheless, we show that the HJB equation employing a stable process as a driving Lévy process admits a smooth and thus classical solution and verify that its associated control policy is optimal. The solution gives useful insights into more complicated cases. In addition, we show a uniqueness result of the HJB equation from a viscosity viewpoint [54] under a wider condition, along with a special treatment of its non-smooth property. More complicated cases are managed numerically.

Our numerical method directly manages the HJB equation without using extra approximations based on the vanishing discount or long-time integration [55-58]. We do not employ these methods but solve the target problem directly. Using a Gauss–Seidel fast-sweep technique [59], the effective Hamiltonian and the potential are iteratively updated until both converge. The non-smooth properties of the HJB equation hinder us from theoretically analyzing convergence of the numerical scheme. Instead, we examine its convergence in terms of both potential and effective Hamiltonian experimentally, focusing on regularity of the non-local term. We also consider an extended model subject to ambiguity [60, 61] numerically. Consequently, we contribute to new mathematical modeling, analysis, and numerical computation of a sustainable environmental engineering problem having unique non-smooth and non-local properties.

The rest of this paper is structured as follows. Our mathematical model is formulated and an exactly-solvable case is analyzed in Section 2. The numerical scheme to discretize the HJB equation is explained in Section 3. HJB equations are numerically computed in Section 4. Summary and future perspectives of this paper are presented in Section 5. Proofs of **Propositions** are placed in **Appendices**.

## 2. Mathematical model
### 2.1 Stochastic process model

We work on a usual complete probability space (e.g., [50]). Key ingredients of the sediment replenishment problem are a sediment storage process as a non-negative and non-increasing continuous-time process and an intervention policy to replenish sediment [44]. See, also **Figure 1**. The time is denoted as $t \geq 0$ and the sediment storage is denoted as $X_t$ at $t$. The range of $X = (X_t)_{t \geq 0}$ is $D = [0,1]$, where $X_t = 0$ represents the depleted state while $X_t = 1$ the fully-replenished state. We assume that the sediment storage decreases by continuous base river flows and/or jump (flood) events [35, 37-38] and that it does not increase unless it is replenished by a decision-maker who manages the environment.

Set a non-decreasing Lévy process [48] $L = (L_t)_{t \geq 0}$ having jumps in $(0, +\infty)$ and the measure $v(dz)$. We consider the infinite activities case $\int_{(0,+\infty)} v(dz) = +\infty$ with $\int_{(0,+\infty)} \min\{1, z\} v(dz) < +\infty$. This assumption is satisfied in major Lévy processes having both small and large jumps such as one-sided stable processes. A state-dependent jump process $\hat{L} = (\hat{L}_t)_{t \geq 0}$ is the jump process having the truncated increment formally expressed as $d\hat{L}_t = \min\{L_t - L_{t-}, X_{t-}\}$. This truncation physically means that the sediment flushed out at each jump event is at most $X_{t-}$.

The sediment storage dynamics without interventions are formulated as

$$dX_t = -S(X_t)\chi_{\{X_t > 0\}} dt - d\hat{L}_t \quad \text{for} \quad t > 0 \quad \text{with} \quad X_0 = x \in [0,1] \tag{1}$$

with a non-negative function $S \in C(D)$. Here, $\chi_A$ is an indicator function of the set $A$: it equals 1 if $A$ is true and 0 otherwise. Sediment is significantly flushed out at floods (large jumps of the second term in the right-hand side of (1)) while gradually at low flows (first term). Here, no clear condition to separate "floods" and "not floods" is imposed. The SDE (1) is therefore a toy model that accounts for both gradual and sudden decrease of the sediment storage. Nevertheless, it is completely new and has a unique non-smoothness due to having the indicator function that ensures non-negativity of $X$ in the absence of jumps. At each jump time $\tau$ of $\hat{L}$, we get the non-negativity

$$X_\tau = X_{\tau-} - (\hat{L}_\tau - \hat{L}_{\tau-}) = X_{\tau-} - \min\{X_{\tau-}, L_\tau - L_{\tau-}\} \geq 0. \tag{2}$$

If one wants to discuss sediment storage dynamics driven purely by jumps, then simply set $S \equiv 0$. We assume that the SDE (1) admits a unique path-wise solution for simplicity.

We assume Poisson opportunities of the discrete and random observation/intervention that the decision-maker can update the information only at jump times of a Poisson process $N = (N_t)_{t \geq 0}$ with an intensity $\Lambda > 0$. The sequence of jump times of $N$ is denoted as $\{\tau_l\}_{l=0,1,2,...}$ where $\tau_0 = 0$ without loss of generality. At each $\tau_l$, the decision-maker can replenish the sediment as $X_{\tau_l+} = X_{\tau_l} + \eta_l$ with some $\eta_l \in \Xi(X_{\tau_l})$, where $\Xi(X_{\tau_l})$ represents the admissible range of $\eta_l$ given $X_{\tau_l}$. In this way, the replenishment is a càglàd (left continuous with right limits) process as in the existing Poisson models [62-64]. We consider the two cases: given $x \in D$, $\Xi_1(x) = \{0, 1-x\}$ (No replenishment or full replenishment) or $\Xi_2(x) = \{\text{If } x = 0 \text{ then } \{0, 1\}. \text{ Otherwise, } \{0\}\}$ (Full replenishment only when the sediment is depleted). In either case, the process $X$ still has the range $D$. Set the sequence of interventions $\bar{\eta} = \{\eta_l\}_{l=0,1,2,\ldots}$ where $\eta_0 = 0$ without loss of generality.

A natural filtration that the collected information at each $\tau_l$ generates is $\mathcal{F} = (\mathcal{F}_t)_{t\geq 0}$ with $\mathcal{F}_t = \sigma\left\{\left(\tau_j, X_{\tau_j}\right)_{0\leq j\leq k}, k = \sup\{j : \tau_j \leq t\}\right\}$. A set of admissible controls, denoted as $C$, contains all sequences $\bar{\eta} = \{\eta_l\}_{l=0,1,2,\ldots}$ such that $\eta_l \in \Xi(X_{\tau_l})$ with $\mathcal{F}_{\tau_l}$-measurable $\eta_l$. Now, the controlled dynamics contain $X_{\tau_l+} = X_{\tau_l} + \eta_l$ at each time $t = \tau_l$ and the SDE (1) at $t \neq \tau_l$. The path-wise uniqueness assumption of $X$ is not violated by this Poisson intervention.

Our control problem is to optimize $\bar{\eta} = \{\eta_l\}_{l=0,1,2,\ldots}$ so that an objective function is minimized. It contains the penalization of the sediment depletion and the cost of replenishment. As we consider a problem in a long-run, a time-averaged objective function is employed:

$$\phi(x, \bar{\eta}) = \liminf_{T\to +\infty} \frac{1}{T} \mathbb{E}^x\left[\int_0^T \chi_{\{X_s=0\}} \mathrm{d}s + \sum_{k\geq 1, \tau_k \leq T}\left(c\eta_k + d\chi_{\{\eta_k>0\}}\right)\right], \qquad (3)$$

where $\mathbb{E}^x$ is the conditional expectation given $X_0 = x \in D$, $c > 0$ and $d > 0$ are the coefficients of proportional and fixed (implementation) costs, respectively. The right-hand side of (3) is bounded since

$$0 \leq \frac{1}{T} \mathbb{E}^x\left[\int_0^T \chi_{\{X_s=0\}} \mathrm{d}s + \sum_{k\geq 1, \tau_k \leq T}\left(c\eta_k + d\chi_{\{\eta_k>0\}}\right)\right] \leq 1 + \frac{1}{T}\cdot\Lambda T(c+d) = 1 + \Lambda(c+d) < +\infty. \qquad (4)$$

Therefore, we can consider a minimization problem of $\phi$ with respect to $\bar{\eta} \in C$. The optimized objective $H : D \to \mathbb{R}$ is set as

$$H(x) = \inf_{\bar{\eta} \in C} \phi(x, \bar{\eta}). \qquad (5)$$

Clearly, $H$ is non-negative and satisfies $H(x) \leq \phi(x, \bar{\eta}_0) = 1$ with a null-control $\bar{\eta}_0 \in C$ such that $\eta_i = 0$ ($i \geq 0$). A minimizing control of (5) is called an optimal control and is denoted as $\bar{\eta} = \eta^*$. The goal of the control problem is to find $\eta^*$.

It is often reasonable to assume that (5) is a non-negative constant following the conventional ergodic control problems [8, 20]. Therefore, we simply write it as $H$. Later, we justify this treatment in **Propositions 1-2**. In the next sub-section, we present an HJB equation that is derived formally. Its justification is also provided by the same propositions.

## 2.2 HJB equation
Applying formal dynamic programming argument [50] leads to the HJB equation

$$H + \mathcal{A}\Phi + \Lambda\left(\Phi - \min_{\eta\in\Xi(x)}\left\{\Phi(x+\eta) + c\eta + d\chi_{\{\eta>0\}}\right\}\right) - \chi_{\{x=0\}} = 0, \quad x \in D, \qquad (6)$$

where the operator $\mathcal{A}$ is defined for generic sufficiently regular $\Psi = \Psi(x)$ as

$$\mathcal{A}\Psi = S(x)\chi_{\{x>0\}}\frac{\mathrm{d}\Psi}{\mathrm{d}x} + \int_{(0,\infty)}\left\{\Psi(x) - \Psi(x - \min\{x, z\})\right\}v(\mathrm{d}z). \qquad (7)$$

We restrict ourselves to functions $\Psi : D \to \mathbb{R}$ such that (7) is well-defined. The function $\Phi : D \to \mathbb{R}$ is called a potential function, which is an auxiliary function to determine $H$ and $\eta^*$ (e.g., see Cao [18]). Using the HJB equation, in the framework of Markov control, we infer an (candidate of) optimal control

$$\eta^*(x) = \arg\min_{\eta\in\Xi(x)}\left\{\Phi(x+\eta) + c\eta + d\chi_{\{\eta>0\}}\right\}, \quad x \in D, \qquad (8)$$

which later turns out to be indeed optimal.

For the HJB equation (6), a solution is a couple $(H, \Phi) \in (\mathbb{R}, C(D))$. A solution is called a classical solution if it satisfies (6) pointwise. Therefore, a classical solution should be $\Phi \in C(D) \cap C^1((0,1])$. If $(H, \Phi)$ is a classical solution, then $(H, \Phi + C)$ with $C \in \mathbb{R}$ also is, implying that an additional condition needs be introduced to compensate this indeterminacy. We impose the following constraint to resolve this issue: we often choose $\tilde{x} = 0$.

$$\Phi(\tilde{x}) = 0 \text{ at one } \tilde{x} \in D. \tag{9}$$

### 2.3 Exact solution: derivation

The HJB equation (6) has the following three non-smooth properties coming from the dynamics and objective function. The non-smooth drift coefficient comes from the non-negativity preserving gradual decrease of the sediment storage, while the non-smooth integrand of (7) from the non-negativity preserving jumps. The last non-smoothness is the source term in (6) on the penalization of the depletion. We show that assuming a stable process as a driving jump process $L$ with $\Xi = \Xi_2$ and $S(x) = \mu x^{1-\alpha}$ ($\mu \geq 0$) leads to an exact classical solution to (6).

Set $v(\mathrm{d}z) = \lambda z^{-(\alpha+1)}\mathrm{d}z$ with $\alpha \in (0,1)$. We thus assume a one-sided stable subordinator having finite variations. The stable process is the simplest Lévy process having infinite activities [48]. Choosing larger $\lambda$ implies stronger river flows, while choosing larger $\alpha$ implies less intermittency. From an engineering viewpoint, choosing larger $\alpha$ models more moderate flows, and vice versa. In particular, the latter characteristic is found in **Figure 2** of controlled sample paths presented later.

For $x \in (0,1]$, we can explicitly calculate the non-local term as

$$\begin{aligned} &\int_{(0,\infty)} \{\Psi(x) - \Psi(x - \min\{x, z\})\} v(\mathrm{d}z) \\ &= \int_{(0,x)} \{\Psi(x) - \Psi(x - z)\} v(\mathrm{d}z) + \int_{(x,\infty)} \{\Psi(x) - \Psi(0)\} v(\mathrm{d}z) \\ &= \int_{(0,x)} \left\{\Psi(x) - \Psi(x - z) - z \frac{\mathrm{d}\Psi(x)}{\mathrm{d}x}\right\} v(\mathrm{d}z) + \frac{\lambda}{1-\alpha} x^{1-\alpha} \frac{\mathrm{d}\Psi(x)}{\mathrm{d}x} + \frac{\lambda}{\alpha} x^{-\alpha} \{\Psi(x) - \Psi(0)\} \end{aligned} \tag{10}$$

This reformulation makes the integral term simpler so that it can be handled more easily in numerical computation; the integral is now considered in a bounded domain and has a compensated form [48]. Instead, new drift and decay terms arise. The coefficient $x^{-\alpha}$ of the decay term blows up as $x \to +0$, but this term will be handled with a simple numerical trick.

Now, we derive an exact classical solution to the HJB equation (6). Set $\Phi(0) = 0$ following (9). This treatment can mitigate the singularity at $x = 0$. The HJB equation in the present case is

$$H + \chi_{\{x>0\}} \mathcal{A}\Phi + \Lambda \left( \Phi - \min_{\eta \in \Xi_2(x)} \{\Phi(x+\eta) + c\eta + d\chi_{\{\eta>0\}}\} \right) - \chi_{\{x=0\}} = 0, \quad x \in D, \tag{11}$$

where

$$\mathcal{A}\Phi = \left(\mu + \frac{\lambda}{1-\alpha}\right) x^{1-\alpha} \frac{\mathrm{d}\Phi(x)}{\mathrm{d}x} + \int_{(0,x)} \left\{\Phi(x) - \Phi(x-z) - z\frac{\mathrm{d}\Phi(x)}{\mathrm{d}x}\right\} \lambda z^{-(\alpha+1)}\mathrm{d}z + \frac{\lambda}{\alpha} x^{-\alpha} \Phi(x). \tag{12}$$

The following proposition is our first mathematical analysis result.

***Proposition 1***

*The HJB equation (11) admits the classical solution $(\hat{H}, \hat{\Phi})$, where*

$$\hat{H} = \frac{1 + \Lambda \min\{c+d, \kappa\}}{1 + \kappa\Lambda} \in \left(\frac{1}{1+\kappa\Lambda}, 1\right] \text{ and } \hat{\Phi}(x) = -\kappa\hat{H}x^\alpha \in C(D) \cap C^1((0,1]) \tag{13}$$

*with*

$$\kappa = \left\{\mu\alpha + \lambda\left(\frac{\alpha}{1-\alpha} + \frac{1}{\alpha} + I_\alpha\right)\right\}^{-1} > 0 \quad \text{and} \quad I_\alpha = \int_{(0,1)} \frac{1-(1-u)^\alpha - \alpha u}{u^{1+\alpha}} du < 0. \qquad (14)$$

**Remark 1** For a compound Poisson case, HJB equations admit solutions discontinuous at $x = 0$ [43, 44]. On the other hand, our explicit solution is continuous. This is considered due to the infinite-activities nature, namely continuous-like nature, of the underlying jumps.

As a by-product of **Proposition 1**, we get

$$\eta^*(X_t) = \begin{cases} 1 & (X_t = 0 \text{ and } c+d \leq \kappa) \\ 0 & (\text{Otherwise}) \end{cases}. \qquad (15)$$

This control is activated if the sediment is depleted and the cost $c+d$ is sufficiently small. **Figure 2** shows sample paths of the controlled sediment storage dynamics for different values of $\alpha$ using (15) and the scheme [65]. A larger $\alpha$ generates paths having more frequent large jumps as mentioned above.

The potential $\hat{\Phi}$ is uniformly bounded in $D$ and is continuously differentiable arbitrary times in $(0,1]$; however, it is only Hölder continuous with the exponent $\alpha$. An important point is that each term of (11), even with an unbounded coefficient, is uniformly bounded in $D$. This finding is a key in verifying optimality of the classical solution and the associated control in the next sub-section.

We close this sub-section with a brief analysis of the (candidates of) effective Hamiltonian $\hat{H}$ and the optimal control $\eta^*$. Clearly, this $\hat{H}$ is a non-deceasing function of the sum of the cost coefficients $c+d$, intuitively meaning that increasing the cost coefficients increases the averaged cost. The excessive increase of the potential cost leads to a null optimal policy where no replenishment is optimal. In addition, $\hat{H}$ is non-increasing with respect to the intensity $\Lambda$, while $\eta^*$ does not explicitly depend on it. Since $\Lambda$ is the inverse of the mean interval of successive observations, this means that observing more frequently can reduce the potential disutility and thus the averaged cost if the decision-maker follow the optimal policy. Notice that the observation cost itself is not included in the model if it is smaller than the implementation cost. If this condition is not satisfied, then the observation cost should be incorporated into the optimization procedure. The obtained $\hat{H}$ also depends monotonically on $\mu, \lambda$ in an increasing manner, explaining that stronger river flows induce more frequent sediment depletion and larger disutility.

### 2.4 Exact solution: verification and comparison results

There are two objectives in this sub-section. The first objective is to verify optimality of the exact solution in **Proposition 1** and its associated control (15). We then show that this control is indeed optimal and that $H = \hat{H}$. The second objective is to show that the classical solution in **Proposition 1** is the unique classical solution with the innocuous constraint $\hat{\Phi}(0) = 0$. We tackle this issue from a viscosity viewpoint [54]. Hence, our uniqueness result is wider than that for classical solutions. In the above-mentioned way, the verification result gives an existence of viscosity solutions, and the comparison result gives their uniqueness. **Proposition 2** shows that the potential and control in **Proposition 1** are optimal.

*Proposition 2*

We have $H = \hat{H}$. In addition, the control (15) is optimal.

**Remark 2** The equation (53) in **Proof of Proposition 2** combined with the uniform boundedness of $\hat{\Phi}$ leads to

$$\hat{H} = \limsup_{T \to +\infty} \frac{1}{T} \mathbb{E}^x \left[ \int_0^T \chi_{\{X_t = 0\}} dt + \sum_{k \geq 1, \tau_k \leq T} \left(c\eta_k^* + d\chi_{\{\eta_k^* > 0\}}\right) \right]. \qquad (16)$$

Therefore, we have the true limit and "liminf" in (3) can be replaced by "lim" or "limsup" when preferred.

***Remark 3*** A similar proof of **Proposition 2** applies to the case with $\Xi = \Xi_1$ if there exists a solution satisfying the corresponding HJB equation in a pointwise sense except at finite number of points in $(0,1]$.

We move to the comparison result. In order to present the result, viscosity solutions to the HJB equation (11) are defined. Before doing that, its non-local term is reformulated so that we can follow [54, Theorem 3]: for any $\varphi \in C^1(D)$ and $x \in (0,1]$, with the notation $0 \le j(x,z) = \min\{x,z\} \le z$ for $(x,z) \in D \times (0, +\infty)$, we have

$$\int_{(0,\infty)} \{\varphi(x) - \varphi(x - \min\{x,z\})\} v(\mathrm{d}z)$$
$$= \int_{(0,\infty)} \left\{\varphi(x) - \varphi(x - j(x,z)) - j(x,z)\frac{\mathrm{d}\varphi(x)}{\mathrm{d}x}\right\} v(\mathrm{d}z) + \int_{(0,\infty)} j(x,z)\frac{\mathrm{d}\varphi(x)}{\mathrm{d}x} v(\mathrm{d}z). \quad (17)$$
$$= \int_{(0,\infty)} \left\{\varphi(x) - \varphi(x - j(x,z)) - j(x,z)\frac{\mathrm{d}\varphi(x)}{\mathrm{d}x}\right\} v(\mathrm{d}z) + \frac{\lambda}{\alpha(1-\alpha)} x^{1-\alpha} \frac{\mathrm{d}\varphi(x)}{\mathrm{d}x}$$

The relationship $0 \le j(x,z) \le z$ means that the last integral of (17) is well-defined.

We define continuous viscosity solutions following the literature of non-local elliptic HJB equations having a singular Lévy measure at $z = 0$ [54, Definition 2 and Proposition 1]. In the definition below, $0 < \delta < 1$ is a constant and $B(x,\delta)$ is an open ball with the center $x \in D$ and the radius $\delta$.

***Definition 1***
***Viscosity sub- and super-solutions***
*A pair of $h \in \mathbb{R}$ and $\psi \in C(D)$ is a viscosity sub-solution (resp., super-solution) if the following condition is satisfied: at each $\hat{x} \in D$, for any $0 < \delta < 1$, for any $\varphi \in C^2(D)$ such that $\psi - \varphi$ is maximized (resp., minimized) at $x = \hat{x}$ on $B(\hat{x},\delta) \cap D$, the following inequality is satisfied:*

$$h + \Lambda\left(\psi(0) - \min_{\eta \in \Xi_2(0)} \{\psi(\eta) + c\eta + d\chi_{\{\eta > 0\}}\}\right) \le 0 \quad (resp.,\ " \ge 0\ ") \text{ if } \hat{x} = 0 \quad (18)$$

*and*

$$h + \left(\mu + \frac{\lambda}{\alpha(1-\alpha)}\right) \hat{x}^{1-\alpha} \frac{\mathrm{d}\varphi(\hat{x})}{\mathrm{d}x}$$
$$+ \int_{(0,\delta)} \left\{\varphi(\hat{x}) - \varphi(\hat{x} - j(\hat{x},z)) - j(\hat{x},z)\frac{\mathrm{d}\varphi(\hat{x})}{\mathrm{d}x}\right\} v(\mathrm{d}z) \quad (resp.,\ " \ge 0\ ") \text{ if } \hat{x} \in (0,1]. \quad (19)$$
$$+ \int_{(\delta,\infty)} \left\{\psi(\hat{x}) - \psi(\hat{x} - j(\hat{x},z)) - j(\hat{x},z)\frac{\mathrm{d}\varphi(\hat{x})}{\mathrm{d}x}\right\} v(\mathrm{d}z) \le 0$$

***Viscosity solutions***
*A pair of $h \in \mathbb{R}$ and $\psi \in C(D)$ is a viscosity solution if it is a viscosity sub-solution as well as a viscosity super-solution.*

***Remark 4*** The condition analogous to Assumption (NLT) [54] is satisfied because we consider the one-sided stable process and $j(x,z) \le z$, meaning that (19) in **Definition 1** is equivalent to conventional definition of viscosity solutions that does not use $\delta$ [54, Proposition 1]. Introducing $\delta$ is therefore not innocuous, and can remove the singularity issue at $z = 0$ owing to the smoothness of the test functions (**Proof of Proposition 3**). Note that the first integral term of (19) is not always well-defined for less regular functions like $\varphi \in C(D)$.

We show a unique solvability result of the HJB equation (11). Our proof follows that of [54,

Theorem 3] and [66, Theorem D.1], but slight modifications are necessary because of its non-smoothness. In addition, we should be careful with scaling the parameters $\varepsilon,\delta$ of the auxiliary function (see, **Appendix A.3**). The proposition below applies not only to the stable case, but also to tempered stable cases with $v(\mathrm{d}z) = \lambda z^{-(\alpha+1)} e^{-bz} \mathrm{d}z$ ($b > 0$).

*Proposition 3*
For any viscosity sub-solution $(H_1, \psi_1)$ and super-solution $(H_2, \psi_2)$, we have $H_1 \leq H_2$. In addition, if $(H_1, \psi_1)$ and $(H_2, \psi_2)$ are viscosity solutions, then $H_1 = H_2$.

By **Propositions 2-3**, we immediately obtain the following theorem.

*Theorem 1*
We have $H = \hat{H}$ as the unique effective Hamiltonian of the HJB equation (11). Furthermore, (15) is an optimal control.

*Remark 5* We have proven the existence (**Proposition 2**) and uniqueness (**Proposition 3**) of viscosity solutions. Of course, the existence itself does not lead to the uniqueness, and vice versa.

## 3. Numerical method
### 3.1 Finite difference scheme
We present a simple iterative numerical method based on the conventional upwind scheme and a fast-sweeping for solving the HJB equation under more generic cases. The scheme discretizes an HJB equation on a uniform grid having the grid points $\{x_i\}_{i=0,1,2,\ldots,M}$ with $M \in \mathbb{N}$ and $x_i = ih$, $h = M^{-1}$.

We explain a discretization process of the HJB equation (6) with $v(\mathrm{d}z) = \lambda z^{-(\alpha+1)} \mathrm{d}z$. Other jumps can be managed similarly. We use (11)-(12) because it has a bounded integration range of the non-local term that is convenient to numerically handle. We assume $\Xi = \Xi_1$ that allows for sediment replenishment for all $x \in D$. The scheme for the simpler case $\Xi = \Xi_2$ is straightforward. Without any loss of generality, set $\Phi(0) = 0$. The approximated $\Phi$ at $x_i$ is denoted as $\Phi_i$.

Each term at $x = x_i$ is discretized as follows. Firstly, $\chi_{\{x=0\}}$ is simply replaced by $\chi_{\{i=0\}}$. The non-local term standing for sediment replenishment is discretized as

$$\Lambda\left(\Phi - \min_{\eta \in \Xi_1(x)}\left\{\Phi(x+\eta) + c\eta + d\chi_{\{\eta>0\}}\right\}\right) \to \Lambda\left(\Phi_i - \min_{k \in \{0, M-i\}}\left\{\Phi_{i+k} + ckh + d\chi_{\{k>0\}}\right\}\right). \tag{20}$$

The optimal control at $x_i$ is then discretized as

$$\eta_i^* = h \times \arg\min_{k \in \{0, M-i\}}\left\{\Phi_{i+k} + ckh + d\chi_{\{k>0\}}\right\}. \tag{21}$$

The drift term is discretized using the upwind scheme as

$$\left(S(x) + \frac{\lambda}{1-\alpha} x^{1-\alpha}\right)\frac{\mathrm{d}\Phi(x)}{\mathrm{d}x} \to \left\{S(x_i) + \frac{\lambda}{1-\alpha} x^{1-\alpha}\right\}\frac{\Phi_i - \Phi_{i-1}}{h}, \tag{22}$$

which is absent if $i = 0$. The decay term arising from jumps is discretized as

$$\frac{\lambda}{\alpha} x^{-\alpha} \Phi(x) \to \frac{\lambda}{\alpha} x_i^{-\alpha} \Phi_i, \tag{23}$$

which is absent if $i = 0$. Finally, the non-local term of jumps is discretized as

$$\int_{(0,x)} \left\{\Phi(x) - \Phi(x-z) - z\frac{\mathrm{d}\Phi(x)}{\mathrm{d}x}\right\} \lambda z^{-(\alpha+1)} \mathrm{d}z \to \sum_{j=1}^{i} I_j \tag{24}$$

with

$$I_j = \left(\Phi_i - \frac{\Phi_{i-j} + \Phi_{i-j+1}}{2} - x_{j+1/2}\frac{\Phi_i - \Phi_{i-1}}{h}\right)\lambda x_{j+1/2}^{-(\alpha+1)}h. \tag{25}$$

Consequently, the discretized HJB equation at the point $x_i$ becomes

$$H + \Lambda\left(\Phi_i - \min_{k \in \{0, M-i\}}\left\{\Phi_{i+k} + ckh + d\chi_{\{k>0\}}\right\}\right) - \chi_{\{i=0\}}$$
$$+ \chi_{\{i>0\}}\left[\left\{S(x_i) + \frac{\lambda}{1-\alpha}x^{1-\alpha}\right\}\frac{\Phi_i - \Phi_{i-1}}{h} + \sum_{j=1}^{i} I_j + \frac{\lambda}{\alpha}x_i^{-\alpha}\Phi_i\right] = 0 \quad \text{for } 0 \leq i \leq M \tag{26}$$

with $\Phi_0 = 0$. We have $M+1$ unknowns ($H$ and $\{\Phi_i\}_{1 \leq i \leq M}$) and $M+1$ equations. The total numbers of the unknowns and the equations are therefore the same.

### 3.2 Fast-sweep method

We apply an iterative method to solve the system of non-linear equations (26). We use a Gauss-Seidel fast-sweep method for Hamilton–Jacobi equations [59] with slight modifications. The main advantages of using a fast-sweep method are its stability against a wide range of Hamilton–Jacobi and related differential equations and its implementability [67, 68] A drawback is that there seems to be no theoretical proof of convergence of the iteration procedure. To the best of the authors' knowledge, this numerical method has not been applied to HJB equations of the ergodic controls like ours.

For $i > 0$, (26) is formally written as

$$H + F_i\Phi_i - G_i\left(\{\Phi_i\}_{1 \leq i \leq M}\right) = 0 \tag{27}$$

with

$$F_i = \Lambda + \chi_{\{i,0\}}\left\{S(x_i) + \frac{\lambda}{1-\alpha}x^{1-\alpha}\right\} + \sum_{j=1}^{i}\lambda z_{j+1/2}^{-(\alpha+1)}h > 0, \tag{28}$$

where $G_i$ represents the terms other than $H + F_i\Phi_i$. The equation for $i = 0$ is

$$H = 1 + \Lambda \min_{k \in \{0, M-i\}}\left\{\Phi_{i+k} + ckh + d\chi_{\{k>0\}}\right\}. \tag{29}$$

Being different from the conventional Hamilton–Jacobi equations that have been solved with fast-sweep methods, our problem has the two qualitatively different unknowns $H$ and $\{\Phi_i\}_{1 \leq i \leq N}$. However, owing to imposing the constraint $\Phi_0 = 0$, we can implement the following pseudo-code to iteratively solve the system (26). The super-script represents iteration steps.

**(Pseudo-code)**
1. Set $n = 0$ and an initial guess $\{\Phi_i^{(0)}\}_{1 \leq i \leq M}$.
2. Set $n \to n+1$.
3. Compute $H^{(n)}$ using (29).
4. Compute $\{\Phi_i^{(n)}\}_{1 \leq i \leq M}$ using $\Phi_i^{(n+1)} = R\Phi_i^{(n)} + (1-R)F_i^{-1}\left\{G_i\left(\{\Phi_i^{(n*)}\}_{1 \leq i \leq M}\right) - H^{(n)}\right\}$ with some relaxation parameter $R \in (0,1)$. The Gauss-Seidel sweep direction is from $i = 1$ to $i = M$. Here, $\Phi_i^{(n*)}$ represents $\Phi_i^{(n+1)}$ if it is available and $\Phi_i^{(n)}$ otherwise.
5. Terminate the iteration if the difference $\left|H^{(n)} - H^{(n-1)}\right|$ and/or $\left\{\left|\Phi_i^{(n)} - \Phi_i^{(n+1)}\right|\right\}_{1 \leq i \leq M}$ are small. Otherwise, go to Step 2.

The relaxation parameter is introduced to enhance computational stability of the iteration because we found that the iteration without relaxation ($R = 0$) diverges.

***Remark 6*** The drift term may be alternatively discretized using a higher-resolution scheme if preferred.

***Remark 7*** Our discretization is not monotone due to the last term of (25). The conventional convergence argument [69] therefore seems not to apply. We therefore check its convergence experimentally.

## 4. Numerical computation
### 4.1 Optimal control and effective Hamiltonian

Because of the lack of theoretical convergence proof, we firstly examine the proposed numerical scheme against the exact solution of **Proposition 1**. Then, the scheme is applied to the more complicated cases. We examine different values of $\alpha$ that critically affects the regularity of the potential. We set the following parameter values in this sub-section: $c = 0.15$, $d = 0.05$, $\mu = 0.10$, $\Lambda = 0.25$, $\lambda = 0.20$. The relaxation parameter is $R = 0.5$. The computational resolution is chosen as $M = 50i$ ($i = 1, 2, 4, 8, 32$) to experimentally analyze convergence of numerical solutions. Iterations of the fast-sweeping are terminated if the maximum point-wise difference between successive updates of the potential becomes smaller than the threshold $10^{-10}$.

**Figures 3** plots the exact and numerical potentials. **Tables 1-3** present the maximum point-wise error of the potential and the absolute error of the effective Hamiltonian for $\alpha$ of 0.2, 0.5, and 0.8, respectively. **Figure 3** shows that the numerical solutions do not have spurious oscillations and qualitatively reproduce the exact solutions. Their convergence is analyzed with **Table 2** more in detail. On the potential $\Phi$, the convergence rates of the numerical solutions are close to 0.2, 0.5, and 0.8 for $\alpha$ of 0.2, 0.5, and 0.8, respectively. This suggests that the convergence rate of the potential in the sense of point-wise maximum is $\alpha$. Convergence rates of $H$ seem to be less affected than that of the potential, but a higher convergence rate is achieved for a higher value of $\alpha$ as in the case of the potential. This parameter dependence is in accordance with the regularity of the potential because smaller $\alpha$ leads to lower regularity and thus slower convergence of numerical schemes. For similar degenerate elliptic equations having non-local terms, dependence of the convergence rates on the regularity of non-local terms has been reported [70, 71]. However, they are different from the presented ones possibly because different equations are discussed in the different works. Nevertheless, the qualitative dependence of the convergence rate on the shape parameter $\alpha$ in our case is natural considering the regularity of $\Phi$.

Hereafter, we set the less restrictive choice $\Xi = \Xi_1$. We apply the presented scheme to numerical computation of the effective Hamiltonian $H$ and the optimal control $\eta^*$ for more complicated than that analyzed in Section 2. We have experimentally found in advance that the optimal control with $\Xi = \Xi_1$ is

$$\eta^*(x) = \begin{cases} 1-x & (0 \leq x \leq \bar{x}) \\ 0 & (\bar{x} < x \leq 1) \end{cases} \qquad (30)$$

with some $\bar{x} \in [0,1)$, or do nothing $\eta^* \equiv 0$. This is an intuitive policy that the sediment should be replenished if its storage is smaller than a threshold. The case $\eta^* \equiv 0$ arises if $c$ and/or $d$ are large. Hence, our focus here is on $H$ and $\bar{x}$.

As a demonstrative example, we analyze parameter dependence of $H$ and $\bar{x}$ on $\alpha$. Before carrying out the numerical analysis, we should normalize $v$ such that different jump processes have a common statistical property so that we can compare their impacts fairer way. Because the jump larger than the size 1 has the same sediment flushing out with that of the size 1, we normalize $v$ so that $\int_1^\infty v(\mathrm{d}z)$ is the same value among $v$ having different $\alpha$. This normalization means that impacts of large jumps are the same among different jump processes. By $\int_1^\infty z^{-(1+\alpha)}\mathrm{d}z = \alpha^{-1}$, we effectively replace $\lambda$ by $\alpha\lambda$. In this way, we get the normalization $\int_1^\infty v(\mathrm{d}z) = 1$. This normalization is only for the computation here.

**Figures 4** shows dependence of the effective Hamiltonian $H$ and the threshold $\bar{x}$ on the parameter $\alpha$. We see that $H$ is increasing with respect to $\alpha$, meaning that the decision-maker experiencing larger $\alpha$ (smaller jumps with higher frequency) yields larger disutility and cost. Even no replenishment is optimal ($H = 1$) for large $\alpha \geq 0.95$. The threshold $\bar{x}$ is unimodal and convex with respect to $\alpha$ and has the maximum near $\alpha = 0.65$. The non-monotone dependence shows that the same threshold value becomes optimal for more than one values of $\alpha$.

### 4.2 Problem with model ambiguity

Any control problems potentially face with model uncertainty such that parameters and/or coefficients are not exactly known. There exist diverse ways of defining the model uncertainty such as the non-linear expectation [72] and ambiguity aversion [73]. A common feature of these methodologies is to consider a worst-case problem where the uncertainty maximizes the objective function to be minimized. This is a zero-sum stochastic differential game between the decision-maker and an opponent serving as the uncertainty, called nature.

We apply the multiplier-robust control formalism that has been proposed for controlling diffusion processes by Hansen and Sargent [60] and later extended to jump-diffusion processes in Anderson et al. [61]. We focus on a case that the jump intensity $\Lambda$ is not exactly known. This case corresponds to the problem where the reference replenishment scheme planned by the decision-maker does not necessarily work as he/she is implemented in a distorted way. The distortion between the reference and actual schemes would arise from communications between the decision-maker and a contractor who transport and supply the sediment. Climatic and environmental factors would also affect implementation of replenishment schemes because of the field-work nature of the problem. The uncertainty issue of the observation intensities has not been discussed so far.

The multiplier-robust control approach replaces $\Lambda$ by $\Lambda a_t$ with a positive process $a = (a_t)_{t \geq 0}$ representing uncertainty. The observation intensity is certain when $a_t = 1$ and uncertain otherwise. The process $a$ is considered as an additional variable to be optimized in the control problem, which is adapted to a natural filtration generated by $(N_t, L_t)_{t \geq 0}$. The objective function (3) is extended to

$$\phi(x, \bar{\eta}, a) = \liminf_{T \to +\infty} \frac{1}{T} \hat{\mathbb{E}}^x \left[ \int_0^T \chi_{\{X_s = 0\}} ds + \sum_{k \geq 1, \tau_k \leq T} \left( c\eta_k + d\chi_{\{\eta_k > 0\}} \right) - \frac{\Lambda}{\gamma} \int_0^T (a_s \ln a_s - a_s + 1) ds \right], \quad (31)$$

where the expectation is now taken with respect to the probability measure generated by $(\hat{N}_t, L_t)_{t \geq 0}$, where $\hat{N}$ is a Poisson process having the intensity $\Lambda a_t > 0$. The last term is a penalization term between the true and believed model since the quantity $a_s \ln a_s - a_s + 1$ is proportional to the relative entropy between them [61]. The parameter $\gamma > 0$ is an ambiguity aversion parameter representing an ambiguity neutral decision-maker $\gamma \to +0$ and a more ambiguity-averse one for larger $\gamma$.

We then consider the game-type optimization problem

$$H = \inf_{\bar{\eta}} \sup_a \phi(x, \bar{\eta}, a). \quad (32)$$

Applying a formal dynamic programming argument analogous to Anderson et al. [61] with the Isaacs condition (the order of "inf" and "sup" can be exchanged in our case) leads to the optimality equation called the Hamilton–Jacobi–Bellman–Isaacs (HJBI) equation:

$$H + \mathcal{A}\Phi - \Lambda \max_{a > 0} \left( -a \left\{ \Phi - \min_{\eta \in \Xi(x)} \left\{ \Phi(x + \eta) + c\eta + d\chi_{\{\eta > 0\}} \right\} \right\} - \frac{a \ln a - a + 1}{\gamma} \right) - \chi_{\{x = 0\}} = 0, \quad x \in D \quad (33)$$

or equivalently

$$H + \mathcal{A}\Phi + \frac{\Lambda}{\gamma} \left( 1 - e^{-\gamma \left( \Phi - \min_{\eta \in \Xi(x)} \{\Phi(x + \eta) + c\eta + d\chi_{\{\eta > 0\}}\} \right)} \right) - \chi_{\{x = 0\}} = 0, \quad x \in D. \quad (34)$$

The worst-case ambiguity $a_t^* = a^*(X_t)$ is

$$a_t^* = e^{-\gamma \left( \Phi(X_t) - \min_{\eta \in \Xi(X_t)} \{\Phi(X_t + \eta) + c\eta + d\chi_{\{\eta > 0\}}\} \right)} = e^{-\gamma \left( \Phi(X_t) - \{\Phi(X_t + \eta^*(X_t)) + c\eta^*(X_t) + d\chi_{\{\eta^*(X_t) > 0\}}\} \right)}. \quad (35)$$

Note that we get the same result when we exchange "max" and "min" in (33) is exchanged; namely, the Isaacs condition is satisfied in the present case. The HJBI equation (33) has a stronger non-linearity than (6), and it reduces to (6) under the ambiguity neutral case $\gamma \to +0$. Again, we set $\Phi(0) = 0$.

We numerically compute the HJBI equation (33) using the proposed numerical method. A difference is that the term on the sediment replenishment is discretized as

$$\frac{\Lambda}{\gamma}\left(1 - e^{-\gamma\left(\Phi - \min_{\eta \in \Xi(x)}\{\Phi(x+\eta) + c\eta + d\chi_{\{\eta>0\}}\}\right)}\right) \to \Lambda\Phi_i - \Lambda\Phi_i + \frac{\Lambda}{\gamma}\left(1 - e^{-\gamma\left(\Phi_i - \min_{k \in [0, N-i]}\{\Phi_{i+k} + ckh + d\chi_{\{k>0\}}\}\right)}\right), \quad x \in D. \quad (36)$$

The redundant representation (36) effectively works to enhance computational stability of the fast-sweep method by utilizing the left-most $\Lambda\Phi_i$ into the term $F_i\Phi_i$.

We examine dependence of the effective Hamiltonian $H$, threshold $\bar{x}$, and the worst-case ambiguity $a^*$ on the ambiguity aversion parameter $\gamma$. The computed potentials are not presented here because they are convex functions like those in the ambiguity-neutral case. In **Figure 5**, we see monotonically increasing dependence of $H$ and $\bar{x}$ on $\gamma$, clearly showing that increasing the ambiguity leads to more frequent sediment replenishment that is possibly more costly. The increasing nature is significant for $\gamma = O(10^0)$ and becomes less significant near $\gamma = 10$. The computational results suggest that, given a set of parameter values, there exists some $\gamma$ around which the optimal policy and the associated effective Hamiltonian sharply change. **Figure 5** also shows that the worst-case observation intensity, which is $\Lambda a^*(x)$, decreases as the ambiguity-aversion nature of the decision-maker becomes heavier. There exists a sharp transition of $a^*(x)$ near the threshold $x = \bar{x}$ for relatively large $\lambda$. This implies that the optimal policy and the worst-case ambiguity is deeply related with each other, and that they should be considered simultaneously.

## 5. Conclusions

We proposed and analyzed an ergodic control problem of a non-smooth and jump-driven dynamical system arising in long-run sediment management. The HJB equation was solved exactly in a simplified case and its optimality was verified. In addition, the solution was characterized as the unique continuous viscosity solution. We also proposed a finite difference scheme for discretizing the HJB equations and experimentally suggested that the scheme is convergent and that its accuracy depends on singularity of the Lévy measure. An extended problem having model ambiguity was analyzed numerically.

A direction of future research is deeper mathematical analysis of the extended problems discussed in Section 4, especially theoretical convergence of numerical solutions in a viscosity sense. A question is whether the conventional convergence argument [69] applies to these problems. It is also interesting to consider the proposed model from a viewpoint of forward-backward SDEs [74] with which an alternative characterization of the optimal control will be found. Such a new characterization can be useful not only for deeper mathematical analysis but also for development of more efficient numerical schemes [75]. However, to find and solve appropriate forward-backward SDEs seems to be a non-trivial task due to the non-smooth dynamics. Coupling with water and biological resources management [11, 76] is also an interesting application. From an engineering viewpoint, other type of jumps will be examined to explore more reasonable descriptions of sediment storage dynamics. The authors are currently tackling a designing issue to find a feasible size ($D$) of sediment storage in a river in Japan by considering more complicated jump processes of the sediment transport based on the proposed model [77]. Existence of invariant measure of the sediment storage, which was implicitly assumed here, should be analyzed in detail under more generic conditions, which is currently undergoing.

We have assumed the unique existence of a (non-negative and square-integrable) solution to the SDE because our aims broadly ranged from theory to simulation. This issue itself is non-trivial because the coefficients do not satisfy the standard continuity condition (e.g., Chapter 4.6 of Capasso and Bakstein [78]). Currently, we are trying to resolve this issue under more generic conditions based on taking a limit of regularized SDEs satisfying the unique existence of a solution. At least in the tempered stable case we obtained an affirmative answer even under model ambiguity [77].

**Appendix A**
*A.1 Proof of Proposition 1*

Firstly, substituting $\Phi = \hat{\Phi}(x) = Bx^\alpha$ with $B \in \mathbb{R}$ into (11) for $x > 0$ yields $B = -\kappa \hat{H}$. Then, the HJB equation at $x = 0$ becomes $\hat{H} = 1 + \Lambda \min\{0, -\kappa \hat{H} + c + d\}$. Rearranging this equation yields the first of (13), and then the second one by $B = -\kappa \hat{H}$. Finally, $\kappa > 0$ follows from

$$\frac{1}{\alpha} + \frac{\alpha}{1-\alpha} + I_\alpha = \int_{(0,+\infty)} \frac{1 - \max\{1-u, 0\}^\alpha}{u^{1+\alpha}} \mathrm{d}u > 0. \tag{37}$$

□

*A.2 Proof of Proposition 2*

The proof follows [50, Theorem 2], but the differences of the dynamics and objective functions must be takin into account. Specifically, we have a uniformly bounded potential, while this assumption is not satisfied in Wang [50]. We have the following two steps to complete the proof. At the first step, we show the inequality $H \geq \hat{H}$ and then $H = \hat{H}$ at the second step.

Choose $x \in D$, $\bar{\eta} \in \mathcal{C}$, and $T > 0$. Our starting point is the Itô's formula:

$$\hat{\Phi}(X_T) - \hat{\Phi}(x) = -\int_0^T S(X_t) \chi_{\{X_t > 0\}} \frac{\mathrm{d}\hat{\Phi}}{\mathrm{d}x}(X_t) \mathrm{d}t \\
+ \sum_{\substack{\hat{L}_t \neq \hat{L}_{t-} \\ 0 \leq t \leq T}} \left(\hat{\Phi}(X_t) - \hat{\Phi}(X_{t-})\right) + \sum_{\substack{N_t \neq N_{t-} \\ 0 \leq t \leq T}} \left(\hat{\Phi}(X_{t+}) - \hat{\Phi}(X_t)\right), \tag{38}$$

where the second and third terms in the right-hand side of (38) represent the summations at all the jump times of the processes $L$ and $N$ in $0 \leq t \leq T$, respectively. Recall that the replenishment is carried out to move the storage $X$ from $X_t$ to $X_{t+}$ at each jump of $N$. The first term is well-defined because of the functional form of $\hat{\Phi}$ (**Proposition 1**). The first summation is understood as [48]

$$\sum_{\substack{\hat{L}_t \neq \hat{L}_{t-} \\ 0 \leq t \leq T}} \left(\hat{\Phi}(X_t) - \hat{\Phi}(X_{t-})\right) = \sum_{\substack{\hat{L}_t \neq \hat{L}_{t-} \\ 0 \leq t \leq T}} \left(\hat{\Phi}(X_t) - \hat{\Phi}(X_{t-}) - (\hat{L}_t - \hat{L}_{t-})\frac{\mathrm{d}\hat{\Phi}}{\mathrm{d}x}(X_{t-})\right) - \int_0^T \frac{\mathrm{d}\hat{\Phi}}{\mathrm{d}x}(X_{t-})\mathrm{d}\hat{L}_t, \tag{39}$$

because the process $L$ has infinitely many small jumps and $X_t - X_{t-} = -(\hat{L}_t - \hat{L}_{t-})$ at its jump times, but is written in the summation form to simplify the description.

We reformulate each term of (39) using the fact that $(\hat{H}, \hat{\Phi})$ is a classical solution. Firstly, we have

$$-\int_0^T S(X_t) \chi_{\{X_t > 0\}} \frac{\mathrm{d}\hat{\Phi}}{\mathrm{d}x}(X_t) \mathrm{d}t \\
= \int_0^T \left\{\Lambda\left(\hat{\Phi}(X_t) - \min_{\eta \in \Xi_2(X_t)}\{\hat{\Phi}(X_t + \eta) + c\eta + d\chi_{\{\eta > 0\}}\}\right) - \chi_{\{X_t = 0\}} + \hat{H}\right\} \mathrm{d}t \tag{40} \\
+ \int_0^T \chi_{\{X_t > 0\}} \int_{(0,+\infty)} \{\hat{\Phi}(X_t) - \hat{\Phi}(\max\{X_t - z, 0\})\} \lambda z^{-(\alpha+1)} \mathrm{d}z \, \mathrm{d}t$$

by (10). In addition, we have

$$\sum_{\substack{N_t \neq N_{t-} \\ 0 \leq t \leq T}} \left(\hat{\Phi}(X_{t+}) - \hat{\Phi}(X_t)\right) = \int_0^T \left(\hat{\Phi}(X_t + \bar{\eta}_t) - \hat{\Phi}(X_t)\right) \mathrm{d}N_t. \tag{41}$$

Here, $\bar{\eta}_t$ equals $\eta_l$ at time $\tau_l$. By (38)-(41), we get

$$\mathbb{E}^x\left[\hat{\Phi}(X_T)\right]-\hat{\Phi}(x)$$
$$=T\hat{H}-\mathbb{E}^x\left[\int_0^T \chi_{\{X_t=0\}}\,dt\right]+\mathbb{E}^x\left[\int_0^T\left\{\Lambda\left(\hat{\Phi}(X_t)-\min_{\eta\in\Xi_2(X_t)}\left\{\hat{\Phi}(X_t+\eta)+c\eta+d\chi_{\{\eta>0\}}\right\}\right)\right\}dt\right]$$
$$+\mathbb{E}^x\left[\int_0^T\int_{(0,+\infty)}\left\{\hat{\Phi}(X_t)-\hat{\Phi}(\max\{X_t-z,0\})\right\}\lambda z^{-(\alpha+1)}\,dz\,dt\right] \quad (42)$$
$$+\mathbb{E}^x\left[\int_0^T\left(\hat{\Phi}(X_t+\bar{\eta}_t)-\hat{\Phi}(X_t)\right)dN_t\right]+\mathbb{E}^x\left[\sum_{L_t\neq L_{t-},0\leq t\leq T}\left(\hat{\Phi}(X_t)-\hat{\Phi}(X_{t-})\right)\right]$$

By the Martingale property of compensated Lévy processes, we get
$$\mathbb{E}^x\left[\sum_{L_t\neq L_{t-},0\leq t\leq T}\left(\hat{\Phi}(X_t)-\hat{\Phi}(X_{t-})\right)\right]=\mathbb{E}^x\left[\int_0^T\int_{(0,+\infty)}\left\{\hat{\Phi}(\max\{X_t-z,0\})-\hat{\Phi}(X_t)\right\}\lambda z^{-(\alpha+1)}\,dz\,dt\right] \quad (43)$$
and
$$\mathbb{E}^x\left[\int_0^T\hat{\Phi}(X_t)\,dN_t\right]=\mathbb{E}^x\left[\int_0^T\Lambda\hat{\Phi}(X_t)\,dt\right]. \quad (44)$$

We then arrive at the equality
$$\mathbb{E}^x\left[\hat{\Phi}(X_T)\right]-\hat{\Phi}(x)$$
$$=T\hat{H}-\mathbb{E}^x\left[\int_0^T\chi_{\{X_t=0\}}\,dt\right]+\mathbb{E}^x\left[\int_0^T\Lambda\left\{\hat{\Phi}(X_t+\bar{\eta}_t)-\min_{\eta\in\Xi_2(X_t)}\left\{\hat{\Phi}(X_t+\eta)+c\eta+d\chi_{\{\eta>0\}}\right\}\right\}dt\right]. \quad (45)$$

Because of
$$\hat{\Phi}(X_t+\bar{\eta}_t)\geq\min_{\eta\in\Xi_2(X_t)}\left\{\hat{\Phi}(X_t+\eta)+c\eta+d\chi_{\{\eta>0\}}\right\}-\left(c\bar{\eta}_t+d\chi_{\{\bar{\eta}_t>0\}}\right), \quad (46)$$

applying (46) to (45) yields
$$\mathbb{E}^x\left[\hat{\Phi}(X_T)\right]-\hat{\Phi}(x)\geq T\hat{H}-\mathbb{E}^x\left[\int_0^T\chi_{\{X_t=0\}}\,dt\right]-\mathbb{E}^x\left[\int_0^T\Lambda\left(c\bar{\eta}_t+d\chi_{\{\bar{\eta}_t>0\}}\right)dt\right]$$
$$=T\hat{H}-\mathbb{E}^x\left[\int_0^T\chi_{\{X_t=0\}}\,dt+\sum_{k\geq 1,\tau_k\leq T}\left(c\eta_k+d\chi_{\{\eta_k>0\}}\right)\right] \quad (47)$$

by the martingale property of the compensated Poisson process. By (47), we get the inequality
$$\hat{H}\leq\frac{1}{T}\mathbb{E}^x\left[\int_0^T\chi_{\{X_t=0\}}\,dt+\sum_{k\geq 1,\tau_k\leq T}\left(c\eta_k+d\chi_{\{\eta_k>0\}}\right)\right]+\frac{\mathbb{E}^x\left[\hat{\Phi}(X_T)\right]-\hat{\Phi}(x)}{T}. \quad (48)$$

Because of the uniform boundedness of $\hat{\Phi}$, we have
$$\hat{H}\leq\liminf_{T\to+\infty}\frac{1}{T}\mathbb{E}^x\left[\int_0^T\chi_{\{X_t=0\}}\,dt+\sum_{k\geq 1,\tau_k\leq T}\left(c\eta_k+d\chi_{\{\eta_k>0\}}\right)\right]=\phi(x;\bar{\eta}). \quad (49)$$

Since $x\in D$ and $\bar{\eta}\in\mathcal{C}$ are arbitrary, we conclude $\hat{H}\leq H$.

We proceed to the second step. Choose the control $\bar{\eta}=\eta^*\in\mathcal{C}$ of (15) and some $x\in D$ and $T>0$. The discussion essentially the same with that in the first step yields
$$\mathbb{E}^x\left[\hat{\Phi}(X_T)\right]-\hat{\Phi}(x)$$
$$=T\hat{H}-\mathbb{E}^x\left[\int_0^T\chi_{\{X_t=0\}}\,dt\right]+\mathbb{E}^x\left[\int_0^T\Lambda\left\{\hat{\Phi}(X_t+\eta_t^*(X_t))-\min_{\eta\in\Xi_2(X_t)}\left\{\hat{\Phi}(X_t+\eta)+c\eta+d\chi_{\{\eta>0\}}\right\}\right\}dt\right], \quad (50)$$

but now
$$\min_{\eta\in\Xi_2(X_t)}\left\{\hat{\Phi}(X_t+\eta)+c\eta+d\chi_{\{\eta>0\}}\right\}=\hat{\Phi}(X_t+\eta_t^*(X_t))+c\eta_t^*(X_t)+d\chi_{\{\eta_t^*(X_t)>0\}}. \quad (51)$$

Therefore, we get
$$\mathbb{E}^x\left[\hat{\Phi}(X_T)\right]-\hat{\Phi}(x)=T\hat{H}-\mathbb{E}^x\left[\int_0^T\chi_{\{X_t=0\}}\,dt+\sum_{k\geq 1,\tau_k\leq T}\left(c\eta_k^*+d\chi_{\{\eta_k^*>0\}}\right)\right] \quad (52)$$

and thus

$$\hat{H} = \frac{1}{T}\mathbb{E}^x\left[\int_0^T \chi_{\{X_t=0\}}\,dt + \sum_{k\geq 1,\tau_k\leq T}\left(c\eta_k^* + d\chi_{\{\eta_k^*>0\}}\right)\right] + \frac{1}{T}\left(\mathbb{E}^x\left[\hat{\Phi}(X_T)\right] - \hat{\Phi}(x)\right). \tag{53}$$

Because the left-hand side is a constant, we get the desired equality

$$\hat{H} = \liminf_{T\to+\infty}\frac{1}{T}\mathbb{E}^x\left[\int_0^T \chi_{\{X_t=0\}}\,dt + \sum_{k\geq 1,\tau_k\leq T}\left(c\eta_k^* + d\chi_{\{\eta_k^*>0\}}\right)\right] = \phi(x;\eta^*) = H. \tag{54}$$

Furthermore, by the first step, we get $\phi(x;\eta^*) \leq \phi(x;\bar{\eta})$ for any $\bar{\eta}\in\mathcal{C}$, showing the optimality.

□

### A.3 Proof of Proposition 3

The second statement immediately follows from the first one. Therefore, we focus on the first statement.

Since the domain $D$ is compact, $\psi_1 - \psi_2$ is maximized at some $\hat{x}\in D$. We should consider the two mutually exclusive cases, which are **Case 1:** $\hat{x}=0$ and **Case 2:** $\hat{x}\in(0,1]$.

**Case 1:** $\hat{x}=0$

By **Definition 1**, we have

$$H_1 + \Lambda\left(\psi_1(0) - \min_{\eta\in\Xi_2(0)}\{\psi_1(\eta) + c\eta + d\chi_{\{\eta>0\}}\}\right) \leq 0 \tag{55}$$

and

$$H_2 + \Lambda\left(\psi_2(0) - \min_{\eta\in\Xi_2(0)}\{\psi_2(\eta) + c\eta + d\chi_{\{\eta>0\}}\}\right) \geq 0. \tag{56}$$

Combining (55) and (56) yields

$$H_1 - H_2 + \Lambda(\psi_1(0)-\psi_2(0)) \leq \Lambda\left(\min_{\eta\in\Xi_2(0)}\{\psi_1(\eta)+c\eta+d\chi_{\{\eta>0\}}\} - \min_{\eta\in\Xi_2(0)}\{\psi_2(\eta)+c\eta+d\chi_{\{\eta>0\}}\}\right). \tag{57}$$

A minimizer of $\min_{\eta\in\Xi_2(0)}\{\psi_i(\eta)+c\eta+d\chi_{\{\eta>0\}}\}$ is denoted as $\eta_i$ ($i=1,2$). Then, we get

$$\begin{aligned}
H_1 - H_2 + \Lambda(\psi_1(0)-\psi_2(0)) &\leq \Lambda\left(\{\psi_1(\eta_1)+c\eta_1+d\chi_{\{\eta_1>0\}}\} - \{\psi_2(\eta_2)+c\eta_2+d\chi_{\{\eta_2>0\}}\}\right)\\
&\leq \Lambda\left(\{\psi_1(\eta_2)+c\eta_2+d\chi_{\{\eta_2>0\}}\} - \{\psi_2(\eta_2)+c\eta_2+d\chi_{\{\eta_2>0\}}\}\right),\\
&= \Lambda(\psi_1(\eta_2)-\psi_2(\eta_2))\\
&\leq \Lambda(\psi_1(0)-\psi_2(0))
\end{aligned} \tag{58}$$

where the last inequality of (58) follows from $\hat{x}=0$. Arranging the above inequality gives the desired result $H_1 - H_2 \leq 0$.

**Case 2:** $\hat{x}\in(0,1]$

By **Case 1**, we can assume

$$\max_{x\in D}\{\psi_1(x)-\psi_2(x)\} > \psi_1(0)-\psi_2(0). \tag{59}$$

As in the standard comparison argument [79], for each $\varepsilon>0$, set

$$\Theta(x,y) = \psi_1(x) - \psi_2(y) - \frac{1}{2\varepsilon}(x-y)^2, \quad x,y\in D. \tag{60}$$

A maximum point of $\Theta$ is denoted as $(x_\varepsilon, y_\varepsilon)$. Again, following the standard technique, we get the estimate under the limit $\varepsilon\to 0$ by choosing a sub-sequence:

$$\frac{1}{2\varepsilon}(x_\varepsilon - y_\varepsilon)^2 \to 0 \quad \text{with}\quad (x_\varepsilon, y_\varepsilon)\to(x_0, x_0) \tag{61}$$

such that $\max_{x\in D}\{\psi_1(x)-\psi_2(x)\} = \psi_1(x_0)-\psi_2(x_0)$ for some $x_0\in D$. In addition, we have

$$\Theta(x_\varepsilon, y_\varepsilon) = \max_{(x,y)\in D\times D}\Theta(x,y) \geq \max_{x\in D}\Theta(x,x) = \max_{x\in D}\{\psi_1(x)-\psi_2(x)\} > \psi_1(0)-\psi_2(0). \tag{62}$$

Therefore, the limit point $\bar{x}$ should be different from 0. This means that, by choosing a sufficiently small $\varepsilon$, there is a constant $\underline{x} \in (0, x_0)$ independent from $\varepsilon$ such that $(x_\varepsilon, y_\varepsilon) \in (\underline{x}, 1] \times (\underline{x}, 1]$. In what follows, we always assume such a sufficiently small $\varepsilon$. Choose a small $\delta \in (0, \underline{x})$.

The rest of the proof is based on the facts that $\psi_1(x) - \left\{ \psi_2(y_\varepsilon) + \frac{1}{2\varepsilon}(x - y_\varepsilon)^2 \right\}$ is maximized at $x = x_\varepsilon$ on $B(x_\varepsilon, \delta) \subset (0,1]$ and that $\psi_2(y) - \left\{ \psi_1(x_\varepsilon) - \frac{1}{2\varepsilon}(x_\varepsilon - y)^2 \right\}$ is minimized at $y = y_\varepsilon$ on $B(y_\varepsilon, \delta) \subset (0,1]$. By **Definition 1**, $\varphi_1(x) = \frac{1}{2\varepsilon}(x - y_\varepsilon)^2$ ($x \in D$) and $\varphi_2(y) = -\frac{1}{2\varepsilon}(x_\varepsilon - y)^2$ ($y \in D$) can be tested against $\psi_1$ at $x = x_\varepsilon$ and $\psi_2$ at $y = y_\varepsilon$, respectively. Notice that

$$\frac{d\varphi_1(x_\varepsilon)}{dx} = \frac{d\varphi_2(y_\varepsilon)}{dy} = \frac{x_\varepsilon - y_\varepsilon}{\varepsilon} \equiv p_\varepsilon \text{ and } \frac{d^2\varphi_1(x_\varepsilon)}{dx^2} = -\frac{d^2\varphi_2(y_\varepsilon)}{dy^2} = \frac{1}{\varepsilon}. \tag{63}$$

Furthermore, for any $x, y$ such that $\underline{x} \leq x \leq y \leq 1$, we have the inequalities

$$\left| x^{1-\alpha} \frac{x-y}{\varepsilon} - y^{1-\alpha} \frac{x-y}{\varepsilon} \right| \leq \frac{1-\alpha}{\min\{x,y\}^\alpha} \frac{|x-y|^2}{\varepsilon} \leq \frac{1-\alpha}{\underline{x}^\alpha} \frac{|x-y|^2}{\varepsilon} \tag{64}$$

and

$$\int_{(0,+\infty)} |j(x,z) - j(y,z)|^2 v(dz) = \int_{(0,x)} |z-z|^2 v(dz) + \int_{(x,y)} |x-z|^2 v(dz) + \int_{(y,+\infty)} |x-y|^2 v(dz)$$

$$\leq \lambda \left\{ \frac{|x-y|}{3\underline{x}^{1+\alpha}} + \frac{1}{\alpha y^\alpha} \right\} |x-y|^2 \tag{65}$$

$$\leq \lambda \underline{C} |x-y|^2$$

with a constant $\underline{C} > 0$ depending on $\underline{x}$ but not on $\delta$. In the second line of (65), we used

$$j(x,z) - j(y,z) = \min\{x,z\} - \min\{y,z\} = 0 \text{ for } 0 < z \leq x \leq y \leq 1. \tag{66}$$

The case $\underline{x} \leq y \leq x \leq 1$ leads to the result.

By the viscosity sub- and super-solution properties, we get

$$H_1 + \left( \mu + \frac{\lambda}{\alpha(1-\alpha)} \right) x_\varepsilon^{1-\alpha} \frac{d\varphi_1(x_\varepsilon)}{dx}$$
$$+ \int_{(0,\delta)} \left\{ \varphi_1(x_\varepsilon) - \varphi_1(x_\varepsilon - j(x_\varepsilon, z)) - j(x_\varepsilon, z) \frac{d\varphi_1(x_\varepsilon)}{dx} \right\} v(dz) \tag{67}$$
$$+ \int_{(\delta,\infty)} \left\{ \psi_1(x_\varepsilon) - \psi_1(x_\varepsilon - j(x_\varepsilon, z)) - j(x_\varepsilon, z) \frac{d\varphi_1(x_\varepsilon)}{dx} \right\} v(dz) \leq 0$$

as well as

$$H_2 + \left( \mu + \frac{\lambda}{\alpha(1-\alpha)} \right) y_\varepsilon^{1-\alpha} \frac{d\varphi_2(y_\varepsilon)}{dy}$$
$$+ \int_{(0,\delta)} \left\{ \varphi_2(y_\varepsilon) - \varphi_2(y_\varepsilon - j(y_\varepsilon, z)) - j(y_\varepsilon, z) \frac{d\varphi_2(y_\varepsilon)}{dy} \right\} v(dz) \tag{68}$$
$$+ \int_{(\delta,\infty)} \left\{ \psi_2(y_\varepsilon) - \psi_2(y_\varepsilon - j(y_\varepsilon, z)) - j(y_\varepsilon, z) \frac{d\varphi_2(y_\varepsilon)}{dy} \right\} v(dz) \geq 0$$

Combining (67) and (68) yields

$$H_1 - H_2 \leq \left(\mu + \frac{\lambda}{\alpha(1-\alpha)}\right)\left\{y_\varepsilon^{1-\alpha}\frac{d\varphi_2(y_\varepsilon)}{dx} - x_\varepsilon^{1-\alpha}\frac{d\varphi_1(x_\varepsilon)}{dy}\right\}$$

$$-\int_{(0,\delta)}\left\{\varphi_1(x_\varepsilon) - \varphi_1(x_\varepsilon - j(x_\varepsilon,z)) - j(x_\varepsilon,z)\frac{d\varphi_1(x_\varepsilon)}{dx}\right\}v(dz)$$

$$+\int_{(0,\delta)}\left\{\varphi_2(y_\varepsilon) - \varphi_2(y_\varepsilon - j(y_\varepsilon,z)) - j(y_\varepsilon,z)\frac{d\varphi_2(y_\varepsilon)}{dy}\right\}v(dz) \quad . \quad (69)$$

$$-\int_{(\delta,\infty)}\left\{\psi_1(x_\varepsilon) - \psi_1(x_\varepsilon - j(x_\varepsilon,z)) - j(x_\varepsilon,z)\frac{d\varphi_1(x_\varepsilon)}{dx}\right\}v(dz)$$

$$+\int_{(\delta,\infty)}\left\{\psi_2(y_\varepsilon) - \psi_2(y_\varepsilon - j(y_\varepsilon,z)) - j(y_\varepsilon,z)\frac{d\varphi_2(y_\varepsilon)}{dy}\right\}v(dz)$$

By (63) and the fact that $\varphi_1, \varphi_2$ are quadratic, we obtain the Taylor expansion

$$\varphi_1(x_\varepsilon) - \varphi_1(x_\varepsilon - j(x_\varepsilon,z)) - j(x_\varepsilon,z)\frac{d\varphi_1(x_\varepsilon)}{dx} = -\frac{1}{2}\{j(x_\varepsilon,z)\}^2\frac{d^2\varphi_1(x_\varepsilon)}{dx^2} = -\frac{1}{2\varepsilon}\{j(x_\varepsilon,z)\}^2 \quad (70)$$

and similarly

$$\varphi_2(y_\varepsilon) - \varphi_2(y_\varepsilon - j(y_\varepsilon,z)) - j(y_\varepsilon,z)\frac{d\varphi_2(y_\varepsilon)}{dy} = -\frac{1}{2}\{j(y_\varepsilon,z)\}^2\frac{d^2\varphi_2(y_\varepsilon)}{dy^2} = \frac{1}{2\varepsilon}\{j(y_\varepsilon,z)\}^2. \quad (71)$$

Then, since $\{j(x_\varepsilon,z)\}^2 + \{j(y_\varepsilon,z)\}^2 \leq 2z^2$, we get

$$-\int_{(0,\delta)}\left\{\varphi_1(x_\varepsilon) - \varphi_1(x_\varepsilon - j(x_\varepsilon,z)) - j(x_\varepsilon,z)\frac{d\varphi_1(x_\varepsilon)}{dx}\right\}v(dz)$$

$$+\int_{(0,\delta)}\left\{\varphi_2(y_\varepsilon) - \varphi_2(y_\varepsilon - j(y_\varepsilon,z)) - j(y_\varepsilon,z)\frac{d\varphi_2(y_\varepsilon)}{dy}\right\}v(dz)$$

$$= \frac{1}{2\varepsilon}\int_{(0,\delta)}\left[\{j(x_\varepsilon,z)\}^2 + \{j(y_\varepsilon,z)\}^2\right]v(dz) \quad (72)$$

$$\leq \frac{\lambda}{\varepsilon}\delta^{2-\alpha}$$

$$\to +0 \text{ as } \varepsilon \to +0$$

with $\delta = \delta(\varepsilon) = \varepsilon^\gamma$ and $\gamma > 0$ sufficiently large, which is subsequently assumed below.

Again by (63), we get

$$-\int_{(\delta,\infty)}\left\{\psi_1(x_\varepsilon) - \psi_1(x_\varepsilon - j(x_\varepsilon,z)) - j(x_\varepsilon,z)\frac{d\varphi_1(x_\varepsilon)}{dx}\right\}v(dz)$$

$$+\int_{(\delta,\infty)}\left\{\psi_2(y_\varepsilon) - \psi_2(y_\varepsilon - j(y_\varepsilon,z)) - j(y_\varepsilon,z)\frac{d\varphi_2(y_\varepsilon)}{dy}\right\}v(dz) \quad . \quad (73)$$

$$= \int_{(\delta,\infty)}\left\{\begin{array}{l}\psi_2(y_\varepsilon) - \psi_1(x_\varepsilon) - \{\psi_2(y_\varepsilon - j(y_\varepsilon,z)) - \psi_1(x_\varepsilon - j(x_\varepsilon,z))\} \\ -\{j(y_\varepsilon,z) - j(x_\varepsilon,z)\}p_\varepsilon\end{array}\right\}v(dz)$$

By the maximizing property of $\Theta$, we have the estimate

$$\Theta(x_\varepsilon - j(x_\varepsilon,z), y_\varepsilon - j(y_\varepsilon,z)) \leq \Theta(x_\varepsilon, y_\varepsilon) \quad (74)$$

and thus

$$\psi_2(y_\varepsilon) - \psi_1(x_\varepsilon) - \{\psi_2(y_\varepsilon - j(y_\varepsilon, z)) - \psi_1(x_\varepsilon - j(x_\varepsilon, z))\}$$
$$\leq \frac{1}{2\varepsilon}\{(x_\varepsilon - j(x_\varepsilon, z)) - (y_\varepsilon - j(y_\varepsilon, z))\}^2 - \frac{1}{2\varepsilon}(x_\varepsilon - y_\varepsilon)^2$$
$$= \frac{1}{2\varepsilon}\{j(x_\varepsilon, z) - j(y_\varepsilon, z)\}\{j(x_\varepsilon, z) - j(y_\varepsilon, z) - 2(x_\varepsilon - y_\varepsilon)\} \quad (75)$$
$$= \frac{1}{2\varepsilon}|j(x_\varepsilon, z) - j(y_\varepsilon, z)|^2 + \{j(y_\varepsilon, z) - j(x_\varepsilon, z)\}p_\varepsilon$$

Substituting (75) into (73) using (65) yields

$$-\int_{(\delta,\infty)}\left\{\psi_1(x_\varepsilon) - \psi_1(x_\varepsilon - j(x_\varepsilon, z)) - j(x_\varepsilon, z)\frac{d\varphi_1(x_\varepsilon)}{dx}\right\}v(dz)$$
$$+\int_{(\delta,\infty)}\left\{\psi_2(y_\varepsilon) - \psi_2(y_\varepsilon - j(y_\varepsilon, z)) - j(y_\varepsilon, z)\frac{d\varphi_2(y_\varepsilon)}{dy}\right\}v(dz)$$
$$\leq \int_{(\delta,\infty)}\left\{\begin{array}{l}\frac{1}{2\varepsilon}|j(x_\varepsilon, z) - j(y_\varepsilon, z)|^2 + \{j(y_\varepsilon, z) - j(x_\varepsilon, z)\}p_\varepsilon \\ -\{j(y_\varepsilon, z) - j(x_\varepsilon, z)\}p_\varepsilon\end{array}\right\}v(dz) \quad (76)$$
$$= \frac{1}{2\varepsilon}\int_{(\delta,\infty)}|j(x_\varepsilon, z) - j(y_\varepsilon, z)|^2 v(dz)$$
$$\leq \frac{1}{2\varepsilon}\int_{(0,\infty)}|j(x_\varepsilon, z) - j(y_\varepsilon, z)|^2 v(dz)$$
$$\leq \frac{\lambda C}{2}\frac{|x_\varepsilon - y_\varepsilon|^2}{\varepsilon}$$

In addition, by (64), as $\varepsilon \to +0$ we get

$$\left|y_\varepsilon^{1-\alpha}\frac{d\varphi_2(y_\varepsilon)}{dx} - x_\varepsilon^{1-\alpha}\frac{d\varphi_1(x_\varepsilon)}{dy}\right| \to 0 \quad (77)$$

by $(x_\varepsilon, y_\varepsilon) \in (\underline{x}, 1] \times (\underline{x}, 1]$ and (64).

Now, by the third line of (72) and (76), from (69) we get the following estimate for sufficiently small $\delta, \varepsilon > 0$:

$$H_1 - H_2 \leq \left(\mu + \frac{\lambda}{\alpha(1-\alpha)}\right)\left\{y_\varepsilon^{1-\alpha}\frac{d\varphi_2(y_\varepsilon)}{dx} - x_\varepsilon^{1-\alpha}\frac{d\varphi_1(x_\varepsilon)}{dy}\right\}$$
$$+\int_{(0,\delta)}\frac{1}{\varepsilon}\left[\{j(x_\varepsilon, z)\}^2 + \{j(y_\varepsilon, z)\}^2\right]v(dz) + \frac{\lambda C}{2}\frac{|x_\varepsilon - y_\varepsilon|^2}{\varepsilon} \quad (78)$$

By (72), using the prescribed $\delta = \delta(\varepsilon) = \varepsilon^\gamma$, taking $\varepsilon \to +0$ obtains the desired inequality $H_1 - H_2 \leq 0$. Combining **Cases 1** and **2** completes the proof.

**Declaration** The authors have no conflict of interests to declare.

**Acknowledgements** JSPS Research Grant 18K01714 and 19H03073, Grant for Environmental Research Projects from the Sumitomo Foundation 203160, and a grant from MLIT Japan for management of seaweed in Lake Shinji support this research.

**Tables**

**Table 1.** Maximum point-wise errors of the potential and the absolute error of the effective Hamiltonian for $\alpha$ of 0.2. Their convergence rates (Conv) are also reported. The exact $H$ is 0.8539356.

| $M$ | 50 | 100 | 200 | 400 | 800 | 1,600 |
|---|---|---|---|---|---|---|
| Error of $\Phi$ | 6.530.E-02 | 5.693.E-02 | 4.960.E-02 | 4.319.E-02 | 3.761.E-02 | 3.275.E-02 |
| Error of $H$ | 6.014.E-04 | 3.154.E-04 | 1.684.E-04 | 9.340.E-05 | 5.540.E-05 | 3.640.E-05 |
| Conv of $\Phi$ | 1.978.E-01 | 1.989.E-01 | 1.994.E-01 | 1.997.E-01 | 1.998.E-01 | |
| Conv of $H$ | 9.312.E-01 | 9.053.E-01 | 8.504.E-01 | 7.536.E-01 | 6.060.E-01 | |

**Table 2.** Maximum point-wise errors of the potential and the absolute error of the effective Hamiltonian for $\alpha$ of 0.5. Their convergence rates are also reported. The exact $H$ is 0.7672192.

| $M$ | 50 | 100 | 200 | 400 | 800 | 1,600 |
|---|---|---|---|---|---|---|
| Error of $\Phi$ | 3.775.E-02 | 2.679.E-02 | 1.898.E-02 | 1.343.E-02 | 9.504.E-03 | 6.722.E-03 |
| Error of $H$ | 1.856.E-03 | 9.788.E-04 | 5.168.E-04 | 2.738.E-04 | 1.478.E-04 | 8.175.E-05 |
| Conv of $\Phi$ | 4.947.E-01 | 4.972.E-01 | 4.985.E-01 | 4.992.E-01 | 4.996.E-01 | |
| Conv of $H$ | 9.230.E-01 | 9.215.E-01 | 9.166.E-01 | 8.897.E-01 | 8.539.E-01 | |

**Table 3.** Maximum point-wise errors of the potential and the absolute error of the effective Hamiltonian for $\alpha$ of 0.8. Their convergence rates are also reported. The exact $H$ is 0.8623599.

| $M$ | 50 | 100 | 200 | 400 | 800 | 1,600 |
|---|---|---|---|---|---|---|
| Error of $\Phi$ | 5.356.E-03 | 3.109.E-03 | 1.796.E-03 | 1.035.E-03 | 5.954.E-04 | 3.423.E-04 |
| Error of $H$ | 1.132.E-03 | 6.041.E-04 | 3.201.E-04 | 1.691.E-04 | 8.913.E-05 | 4.713.E-05 |
| Conv of $\Phi$ | 7.846.E-01 | 7.917.E-01 | 7.955.E-01 | 7.972.E-01 | 7.985.E-01 | |
| Conv of $H$ | 9.061.E-01 | 9.162.E-01 | 9.205.E-01 | 9.241.E-01 | 9.192.E-01 | |

**Figures**

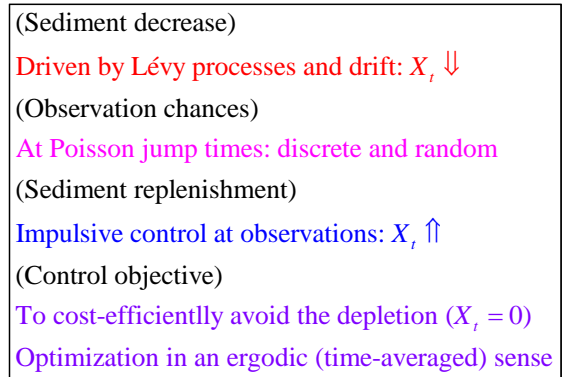

**Figure 1.** Summary of the proposed model.

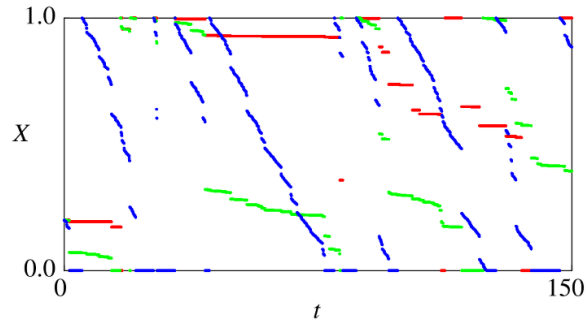

**Figure 2.** Sample paths of the controlled sediment storage dynamics: $\alpha = 0.2$ (Red), 0.5 (Green), and 0.8 (Blue). We used the time increment 0.01 and $X_0 = 0.5$, $\Lambda = 0.15$, and $\lambda = 0.01$.

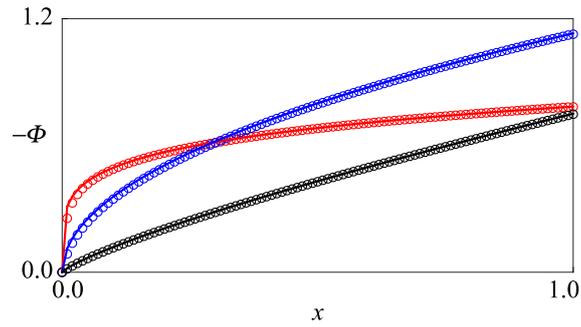

**Figure 3.** Exact (Lines) and numerical potentials (Circles) for $\alpha$ of 0.2 (Red), 0.5 (Blue), and 0.8 (Black). The resolution is $M = 100$.

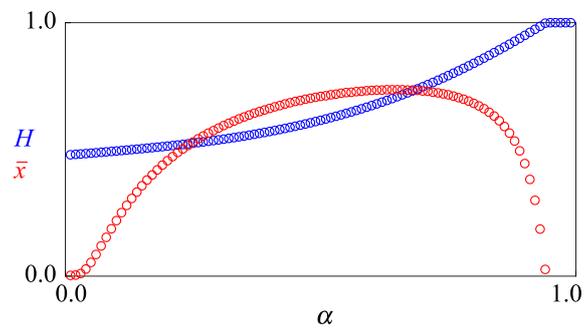

**Figure 4.** The effective Hamiltonian $H$ (Blue) and the threshold $\bar{x}$ (Red) for $0.01 \leq \alpha \leq 0.99$. Notice that $\bar{x}$ does not exists for $\alpha \geq 0.95$ where no replenishment is optimal. The resolution is $M = 800$.

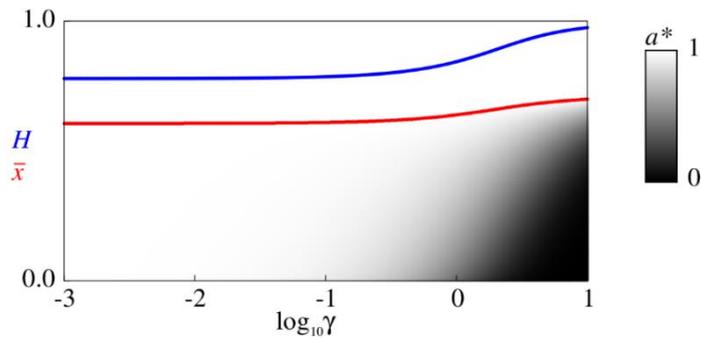

**Figure 5.** The effective Hamiltonian $H$ (Blue), threshold $\bar{x}$ (Red), and The worst-case ambiguity $a = a^*(x)$ (Grey plot) for $10^{-3} \leq \gamma \leq 10^1$. The resolution is $M = 800$.